\documentclass{article}
\usepackage{mathrsfs}
\usepackage{amsfonts}
\usepackage{amsfonts,amssymb,amsmath}

\makeatletter
\@addtoreset{equation}{section}

\makeatother
\newtheorem{thm}{\hskip 20pt Therorem}[section]

\newtheorem{prop}{\hskip 20pt Proposition}[section]

\newcommand{\R}{\mathbb{R}}

\begin{document}
\title{Different Asymptotic Spreading Speeds Induced by Advection in a Diffusion Problem
with Free Boundaries
\thanks{This work is supported by NSFC (11271285, 11071209).} }
\author{Hong Gu\\
{\small\it Department of Mathematics, Tongji University, Shanghai 200092, China}\\
Zhigui Lin\\
{\small\it School of Mathematical Science, Yangzhou University, Yangzhou 225002, China}\\
Bendong Lou\thanks{Corresponding author. Email: blou@tongji.edu.cn}\\
{\small\it Department of Mathematics, Tongji University, Shanghai 200092, China}}

\date{}
\maketitle
\baselineskip 18pt

\begin{abstract}
In this paper, we consider a Fisher-KPP equation with an advection term and two free boundaries,
which models the behavior of an invasive species in one dimension space.
When spreading happens (that is, the solution converges to a positive constant), we
use phase plane analysis and upper/lower solutions to prove that the rightward and leftward
asymptotic spreading speeds exist, both are positive constants. Moreover, one of them is bigger
and the other is smaller than the spreading speed in the corresponding problem
without advection term.
\end{abstract}

\section{Introduction}

In 2010, Du and Lin  \cite{DuLin} studied the following Fisher-KPP problem with free boundaries:
\begin{equation}\label{p0}
\left\{
\begin{array}{ll}
 u_t- d u_{xx}=u(1-u), &  g(t)< x<h(t),\ t>0,\\
 u(t,g(t))=0,\ \ g'(t)=-\mu u_x(t, g(t)), & t>0,\\
 u(t,h(t))=0,\ \ h'(t)=-\mu u_x (t, h(t)) , & t>0,\\
-g(0)=h(0)= h_0,\ \ u(0,x) =u_0 (x),& -h_0\leq x \leq h_0,
\end{array}
\right.
\tag{$P_{0}$}
\end{equation}
where $d$ and $\mu$ are positive constants, the initial function $u_{0}(x)$ satisfies
\begin{equation}\label{u0}
    u_{0}\in C^{2}([-h_{0},h_{0}]),\ \ u_{0}(\pm h_{0})=0\ \ \textrm{and} \ \ u_{0}>0 \ \
    \textrm{in} \ \ (-h_{0},h_{0}),
\end{equation}
for some $h_0>0$. They used \eqref{p0} to model the spreading of a new or invasive species
with population density $u(t,x)$ over a one dimensional habitat, with the free boundaries
$x=g(t),\ h(t)$ representing the expanding fronts. They obtained a dichotomy result,
that is, either spreading happens ($u(t,\cdot)\to 1$ locally uniformly in $\R$ and
$h(t),\; -g(t)\to \infty$ as $t\to \infty$) or vanishing happens ($u(t,\cdot)\to 0$ uniformly
in $[g(t),h(t)]$ as $t\to \infty$ and $h(t)-g(t)<\infty$).
Furthermore, when spreading happens, they obtained the existence of the asymptotic spreading speed
(\cite[Proposition 4.1]{DuLin}):
\begin{equation}\label{speed-000}
c^*:= \lim_{t\to\infty}\frac{h(t)}{t}  = \lim_{t\to\infty}\frac{-g(t)}{t} >0.
\end{equation}
Recently, further extensions have been done, for example, Du and Guo \cite{DuGuo, DuGuo2}
studied the problem in higher dimension spaces and in heterogeneous environment.
Du and Lou \cite{DuLou} studied the problem with general nonlinear $f$, including
general monostable, bistable and combustion types of $f$. Among others,
they all proved that the asymptotic spreading speed when spreading happens is the same positive constant in any direction.

However, some species prefers to move towards one direction because of rich resource, appropriate climate, etc.
Some diseases spread along the wind direction. In 2009, Maidana and Yang in \cite{Maiyang} studied the propagation of West
Nile Virus from New York City to California state. It was observed that West Nile Virus appeared for the first time in
New York City in the summer of 1999. In the second year the wave front travels 187km to the north and 1100km to the south.
Therefore, they took account of the advection movement and showed that bird advection becomes an important factor for lower mosquito biting rates.
Recently, Averill in \cite{Ave} considered the effect of intermediate advection on the dynamics of two-species competition system, and provides a concrete range of advection strength for the coexistence of two competing species. Moreover, three different kinds of transitions from small advection to large advection were illustrates theoretically and numerically.

What is the difference between the asymptotic spreading speed of the left frontier and that of the right frontier when invasive species is spreading?
To address the question, in this paper we study the following problem with an advection term:
\begin{equation}\label{p1}
\left\{
\begin{array}{ll}
 u_t-du_{xx}+\beta u_{x}=u(1-u), &  g(t)< x<h(t),\ t>0,\\
 u(t,g(t))=0,\ \ g'(t)=-\mu u_x(t, g(t)), & t>0,\\
 u(t,h(t))=0,\ \ h'(t)=-\mu u_x (t, h(t)) , & t>0,\\
-g(0)=h(0)= h_0,\ \ u(0,x) =u_0 (x),& -h_0\leq x \leq h_0,
\end{array}
\right.
\tag{$P_{1}$}
\end{equation}
where $d, \mu, h_0, u_0$ are as above and $\beta >0$ is a constant.

By a similar argument as in \cite{DuGuo, DuLin, DuLou}, we have the following basic results.
\begin{itemize}
 \item[{\rm (i)}] Problem \eqref{p1} has a time global solution $(u,g,h)$ with
$u\in C^{1+\alpha/2,2+\alpha}((0,\infty)$ $\times[g(t),h(t)])$ and $g,h\in C^{1+\alpha/2}([0,\infty))$
for any $\alpha\in(0,1)$;
  \item[{\rm (ii)}] $0<u(t,x)\leq C_1$ for $g(t)<x<h(t),\ t>0$ and $0<-g'(t),h'(t)<C_2$ for $t>0$,
where $C_1$ and $C_2$ are constants independent of $t$.
\end{itemize}
In a forthcoming paper \cite{GLL}, we studied the asymptotic behavior of
the solutions of \eqref{p1}. More precisely, we gave some sufficient conditions for spreading
and some sufficient conditions for vanishing.
It turns out that spreading happens only if
\begin{equation}\label{cond}
0<\beta<2\sqrt{d}.
\end{equation}


This paper is devoted to the difference between the leftward and rightward asymptotic
spreading speeds induced by the advection term $\beta u_x$.
This is an interesting problem from ecological point of view.


\begin{thm}\label{speed}
Assume $0<\beta<2\sqrt{d}$. Let $(u,g,h)$ be a solution of \eqref{p1} for which
spreading happens.  Then the leftward and rightward asymptotic spreading speeds exist:
$$
c_{l}^{*} := \lim_{t\to\infty}\frac{- g(t)}{t},\ \qquad c_r^* := \lim_{t\to\infty}\frac{h(t)}{t}.
$$
Moreover, $0 < c_l^* <c^* <c_r^* $, where $c^{*}$ is the spreading speed of the solution of
\eqref{p0}.
\end{thm}

\noindent
$c^*$ is given in \eqref{speed-000} which is nothing but $k_{0}$ in \cite[Proposition 4.1]{DuLin},
or $c^*$ in \cite[Theorem 1.10]{DuLou}. It depends on the parameter $\mu$.
Similarly $c_{l}^{*}$ and $c_{r}^{*}$ depend on $\mu$. Clearly, $c^*_r$ and $c^*_l$
also depend on $\beta$. On these dependence we have the following results.


\begin{thm}\label{relationship}
\begin{enumerate}
  \item[{\rm (i)}] If $\beta\in (0,2\sqrt{d})$ is fixed, then $c_{l}^{*},c^{*},c_{r}^{*}$ are
  strictly increasing in $\mu$, and
  \begin{eqnarray*}
    \lim_{\mu\to0}c_{l}^{*}=0, & & \lim_{\mu\to\infty}c_{l}^{*}=2\sqrt{d}-\beta, \\
    \lim_{\mu\to0}c^{*}=0, & & \lim_{\mu\to\infty}c^{*}=2\sqrt{d},\\
    \lim_{\mu\to0}c_{r}^{*}=0, & & \lim_{\mu\to\infty}c_{r}^{*}=2\sqrt{d}+\beta;
  \end{eqnarray*}
  \item[{\rm (ii)}] if $\mu$ is fixed, then $c^*_r,\ -c^*_l$ are strictly increasing in $\beta$, and
  $$
  \lim_{\beta\to0}c_{l}^{*}=\lim_{\beta\to0}c_{r}^{*}=c^{*},\ \ \lim_{\beta\to2\sqrt{d}}c_{l}^{*}=0.
  $$
\end{enumerate}
\end{thm}


\section{Semi-waves and spreading speeds}\label{sec:speed}

Throughout this section we assume that \eqref{cond} holds and that $(u,g,h)$ is a solution of
\eqref{p1} for which spreading happens, that is, $h(t), \; -g(t)\to \infty \ (t\to \infty)$,
and $u(t,\cdot)\rightarrow 1$ locally uniformly in $\mathbb{R}$. Denote $f(u):=u(1-u)$
for convenience. We remark that the approaches below remain valid for general
monostable nonlinear $f$.

To determine the spreading speed, we will
construct upper and lower solutions based on semi-waves.

\subsection{Phase plane analysis and semi-waves}\label{ss:semi-waves}

We call $q(z)$ a semi-wave with speed $c$ if $(c,q(z))$ satisfies
\begin{equation}\label{prob-q-infty}
\left\{
\begin{array}{l}
q'' -\frac{c-\beta}{d}q' +\frac{f(q)}{d}=0 \quad \mbox{ for } z\in (0,\infty),\\
q(0)=0, \ q(\infty)=1, \ q(z)>0\ \mbox{ for } z\in (0,\infty).
\end{array}
\right.
\end{equation}
The first equation in this problem is equivalent to the following system:
\begin{equation}\label{q-p}
\left\{
\begin{array}{l}
q'=p,\\ p'=\frac{c-\beta}{d}p-\frac{f(q)}{d}.
\end{array}
\right.
\end{equation}
A solution $(q(z), p(z))$ of this system traces out a
trajectory in the $q, p$-plane or, as it is usually called, the phase plane (cf.
\cite{AW1, AW2, DuLou, Pet}). Such a trajectory has slope
\begin{equation}\label{P}
\frac{\mathrm{d}p}{\mathrm{d}q}=\frac{c-\beta}{d}-\frac{f(q)}{dp}
\end{equation}
at any point where $p\neq0$. Here we are only interested in a trajectory of \eqref{q-p} that
starts from the point $(0,\omega)$ with some $\omega \geq0$ and ends at the
point $(1,0)$ as $z\rightarrow+\infty$.

For any fixed $c\geq0$, $(0,0)$ and $(1,0)$ are critical points of the system \eqref{q-p}.
The eigenvalues of the corresponding linearizations are
$$
\lambda_{0}^{\pm}=\frac{c-\beta\pm\sqrt{(c-\beta)^{2}-4d}}{2d}\ (\mbox{at } (0,0)),
$$
$$
\lambda_{1}^{\pm}=\frac{c-\beta\pm\sqrt{(c-\beta)^{2}+4d}}{2d}\ (\mbox{at } (1,0)),
$$
respectively. Thus $(1,0)$ is a saddle point and $(0,0)$ is
\begin{enumerate}
  \item[(i)] a center or a spiral point, if $0\leq c<\beta+2\sqrt{d}$;
  \item[(ii)] a nodal point, if $c\geq\beta+2\sqrt{d}$.
\end{enumerate}
Therefore, by the theory of ODE (cf. \cite{Pet}), there exactly two trajectories of \eqref{q-p}
that approach $(1,0)$ from $q<1$. One of them, denoted by $T^{c}_{r}$, has slope $\lambda^-_1 <0$
at $(1,0)$. Suppose that $T^{c}_{r}$ is expressed by a function $p=P^{c}_{r}(q)$. Then
$p=P^{c}_{r}(q)$ satisfies \eqref{P} and $T^{c}_{r}$ lies in the semistrip
$$
S=\{(q,p):0<q<1,p>0\}.
$$
$T^{c}_{r}$ is a trajectory through $(1,0)$ and $(0,P^c_r (0^+))$ for some $P^c_r (0^+)\geq0$.
The following are well known results (cf. \cite{AW1, AW2, DuLou, VVV}).

\begin{prop}\label{prop:right-tw}
Let $c_r^0 := 2\sqrt{d} +\beta$. Then

\begin{itemize}
 \item[(i)] for any $c\in [0,c^0_r)$, $P^c_r (0)$ is positive, continuous, strictly
 decreasing in $c$, and $\lim\limits_{c \nearrow c^0_r } P^c_r (0)=0$;

 \item[(ii)] for any $c\geq c^0_r$, $P^c_r (0^+)=0$.
 \end{itemize}
\end{prop}

In case (ii), each $T^c_r$ is a trajectory in $S$ through $(0,0)$ and $(1,0)$, and so
it corresponds to a traveling wave with speed $c$, $c_{r}^{0}$ is nothing but the
minimal speed of these traveling waves.

Denote $\zeta(c):=P^{c}_{r}(0)-\frac{c}{\mu}$ for $c\in [0,c^0_r)$. In view of Proposition \ref{prop:right-tw},
$\zeta(c)$ is continuous and strictly decreasing in $c\in [0,c^0_r)$, and it satisfies
\begin{eqnarray*}
  &&\zeta(0)=P^0_{r} (0)>P^\beta_{r} (0) = \sqrt{\frac{2}{d} \int_0^1 f(s)ds} >0, \\
  &&\zeta((c^0_r)^-)=-\frac{c^0_r}{\mu} <0.
\end{eqnarray*}
Thus there exists a unique $c^{*}_{r}\in(0,c^0_r)$ such that $\zeta(c^{*}_{r})=0$, i.e.
$P_{r}^{c^{*}_{r}}(0)=c^{*}_{r}/\mu$.
Summarizing the above results we have the following proposition.

\begin{prop}\label{cr*}
Problem \eqref{prob-q-infty} has exactly one solution
$(c,q) = (c^{*}_{r},q_{r}^{*})$ such that
\begin{equation}\label{c-q-r}
\mu (q_r^*)'(0)=c^*_r .
\end{equation}
Moreover, $c^{*}_{r}\in (0,\beta+2\sqrt{d})$.
\end{prop}

Later we will use this semi-wave to estimate the rightward spreading speed.
Similarly, to estimate the leftward spreading speed, we will need another semi-wave
traveling to left, which is a solution of the following problem:
\begin{equation}\label{prob-q-infty-l}
\left\{
\begin{array}{l}
q'' -\frac{c+\beta}{d}q' +\frac{f(q)}{d}=0 \quad \mbox{ for } z\in (0,\infty),\\
q(0)=0, \ q(\infty)=1, \ q(z)>0\ \mbox{ for } z\in (0,\infty).
\end{array}
\right.
\end{equation}
Similar as above, this problem can be studied by considering the problem
\begin{equation}\label{P-l}
\frac{\mathrm{d}p}{\mathrm{d}q}=\frac{c+\beta}{d}-\frac{f(q)}{dp}
\end{equation}
in the $q, p$-phase plane, where $p=q'$.
Denote $P^{c}_{l}(q)$ the solution of this equation whose trajectory
through $(1,0)$ and $(0,P^c_l(0^+))$ for some $P^c_l(0^+)\geq 0$.
In a similar way as above we have the following results.

\begin{prop}\label{prop:left-tw}
Let $c_l^0 := 2\sqrt{d} -\beta$. Then

\begin{itemize}
 \item[(i)] for any $c\in [0,c^0_l)$, $P^c_l (0)$ is positive, continuous, strictly
 decreasing in $c$, and $\lim\limits_{c \nearrow c^0_l} P^c_l (0)=0$;

 \item[(ii)] for any $c\geq c^0_l$, $P^c_l (0^+)=0$.
 \end{itemize}
\end{prop}

\begin{prop}\label{cl*}
Problem \eqref{prob-q-infty-l} has exactly one solution
$(c,q) = (c^{*}_l,q_l^{*})$ such that
\begin{equation}\label{c-q-l}
\mu (q_l^*)'(0)=c^*_l .
\end{equation}
Moreover, $c^{*}_l\in(0,c^0_l)$.
\end{prop}

Next, we make suitable perturbations of $f(u)$ to derive corresponding semi-waves that
can be used to construct upper and lower solutions of \eqref{p1}.

For any small $\varepsilon>0$, set
$$
\widetilde{f}_{\varepsilon} (u) := f(u)- \frac{\varepsilon}{1-\varepsilon}u^2 \equiv  u\Big(
1- \frac{1}{1-\varepsilon}u\Big),
$$
$$
\widehat{f}_{\varepsilon} (u) := f(u) + \frac{\varepsilon}{1+\varepsilon}u^2 \equiv  u\Big(
1- \frac{1}{1+\varepsilon}u\Big).
$$
Note that $\widetilde{f}_{\varepsilon} (u)$ is strictly decreasing in $\varepsilon$ and it
has exactly two zeros $0$ and $1-\varepsilon$.
$\widehat{f}_{\varepsilon} (u)$ is strictly increasing in $\varepsilon$ and it
has exactly two zeros $0$ and $1+\varepsilon$.
In a similar way as above, we know that problem \eqref{prob-q-infty} with $f$
replaced by $\widetilde{f}_{\varepsilon}$ (resp. $\widehat{f}_\varepsilon$)
has exactly one solution $(\widetilde{c}_{r}^{*},\widetilde{q}_{r}^{*})$ with
$\mu (\widetilde{q}_r^*)'(0)= \widetilde{c}^*_r$ and $\widetilde{c}^{*}_{r}\in(0,c^0_r)$
(resp. a solution $(\widehat{c}_r^{*},\widehat{q}_r^{*})$ with
$\mu (\widehat{q}_r^*)'(0)= \widehat{c}^*_r$ and $\widehat{c}^{*}_{r}\in(0,c^0_r)$),
where $c_{r}^{0}=2\sqrt{d}+\beta$.
Similarly, problem \eqref{prob-q-infty-l} with $f$ replaced by
$\widetilde{f}_{\varepsilon}$ (resp. $\widehat{f}_\varepsilon$)
has exactly one solution $(\widetilde{c}_l^{*},\widetilde{q}_l^{*})$ with
$\mu (\widetilde{q}_l^*)'(0)= \widetilde{c}^*_l$ and $\widetilde{c}^{*}_l\in(0,c^0_l)$
(resp. a solution $(\widehat{c}_l^{*},\widehat{q}_l^{*})$ with
$\mu (\widehat{q}_l^*)'(0)= \widehat{c}^*_l$ and $\widehat{c}^{*}_l\in(0,c^0_l)$),
where $c_{l}^0=2\sqrt{d}-\beta$.


\begin{prop}\label{cl-cr}
The following conclusions hold.
$$\widetilde{c}_{r}^{*}<c_{r}^{*}<\widehat{c}_{r}^{*},\quad
\lim_{\varepsilon\to0}\widetilde{c}_{r}^{*}=\lim_{\varepsilon\to0}\widehat{c}_{r}^{*}=c_{r}^{*},$$
and
$$\widetilde{c}_{l}^{*}<c_{l}^{*}<\widehat{c}_{l}^{*},\quad
\lim_{\varepsilon\to0}\widetilde{c}_{l}^{*}=\lim_{\varepsilon\to0}\widehat{c}_{l}^{*}=c_{l}^{*}.$$
\end{prop}

{\it Proof}.\ \
We first prove $\widetilde{c}_{r}^{*}<c_{r}^{*}$.
For any $c\in [0, c^0_r)$, consider the problem \eqref{P} with $f$ replaced by $\widetilde{f}_\varepsilon$,
denote the solution with trajectory through the critical point $(0, 1-\varepsilon)$ in the
phase plane by $\widetilde{P}^c_{r,\varepsilon} (q)$.
Similar as Proposition \ref{prop:right-tw} (i) we have $\widetilde{P}^c_{r,\varepsilon} (0)>0$ for all
$c\in [0,c^0_r)$. Moreover, $\widetilde{P}^{c}_{r,\varepsilon}(q)<P^{c}_{r}(q)\ (q\in (0,1-\varepsilon])$
by $\widetilde{f}_\varepsilon(q)\leq f(q)\ (0< q \leq 1-\varepsilon)$.
We now prove
\begin{equation}\label{tildeP<P}
 0< \widetilde{P}^{c}_{r,\varepsilon}(0) <P^{c}_{r}(0)\quad \mbox{for } c\in [0,c^0_r).
\end{equation}
Otherwise, $\widetilde{P}^{c}_{r,\varepsilon}(0) =P^{c}_{r}(0)$, and so the function
$\eta(q):= P^{c}_{r}(q) - \widetilde{P}^{c}_{r,\varepsilon}(q)$ satisfies
$$
\eta' <a(q) \eta\ \ (0<q<1-\varepsilon),\qquad \eta(0)=0,
$$
where $a(q):= \widetilde{f}_\varepsilon (q) [d P^{c}_{r}(q) \widetilde{P}^{c}_{r,\varepsilon}(q)]^{-1}$.
This implies that $\eta(q)<0 \ (0<q<1-\varepsilon)$, a contradiction.

Denote $\tilde{\zeta}(c):=\widetilde{P}^{c}_{r,\varepsilon}(0)-\frac{c}{\mu}$. Then \eqref{tildeP<P} implies that
$$
\tilde{\zeta} (c)< \zeta(c) \quad \mbox{for } c\in [0,c^0_r).
$$
Similar as above, both $\tilde{\zeta}(c)$ and $\zeta(c)$ are continuous and strictly decreasing
functions in $[0,c^0_r)$, and
$$
\tilde{\zeta}((c^0_r)^-)=\zeta((c^0_r)^-)= - \frac{c^0_r}{\mu}.
$$
Therefore $\widetilde{c}_{r}^{*}<c_{r}^{*}$ by their definitions: $\tilde{\zeta}
(\widetilde{c}^*_r)=\zeta(c^*_r)=0$.

Next we prove $\lim_{\varepsilon\to0}\widetilde{c}_{r}^{*}=c_{r}^{*}$. It is sufficient to
show that, for any $c\in[0,c_r^0)$,
\begin{equation}\label{epsilon to 0}
\widetilde{P}^{c}_{r,\varepsilon}(0)\rightarrow P^{c}_{r}(0) \mbox{ as }
\varepsilon\rightarrow0.
\end{equation}
By the monotonicity of $\widetilde{f}_\varepsilon$, it is easily seen that
$\widetilde{P}^{c}_{r,\varepsilon}(q)$ is monotonically decreasing in $\varepsilon$,
and it is bounded from above by $P^{c}_{r}(q)$. Therefore, as $\varepsilon \to 0$,
$\widetilde{P}^{c}_{r,\varepsilon}(q)$ converges to some function $R(q)$
in $C^1([0,1-\delta])$ for any $0<\delta <1$.
Clearly, $p=R(q)$ corresponds to a trajectory of \eqref{q-p} that approaches $(1,0)$ in the phase plane
with a non-positive slope at $(1,0)$. Consequently, $R(q) \equiv P^{c}_{r}(q)$, and so
\eqref{epsilon to 0} is proved.

Other conclusions can be proved in a similar way as above.
\hfill $\square$

\subsection{Asymptotic spreading speed}

 {\it Proof of Theorem \ref{speed}}. We first estimate the rightward asymptotic spreading speed.
For any small $\varepsilon>0$ we define
$$
\widetilde{w} (t,x):=\widetilde{q}_{r}^{*}(\widetilde{c}_{r}^{*}t -x),\ \ x\in[0,\widetilde{c}^*_r t],
$$
Since $(\widetilde{c}^{*},\widetilde{q}_{r}^{*}(z))$ satisfies
\begin{equation*}
\left\{
\begin{array}{l}
q'' -\frac{c-\beta}{d}q' +\frac{\widetilde{f}_{\varepsilon}(q)}{d}=0 \quad \mbox{ for } z\in (0,\infty),\\
q(0)=0, \ q(\infty)=1-\varepsilon, \ q'(0)=\frac{c}{\mu} \mbox{ and } q'(z)>0\ (z>0),
\end{array}
\right.
\end{equation*}
we have
$$
\widetilde{w} (t,x)\leq1-\varepsilon,\ \quad
\widetilde{w}_t- \widetilde{w}_{xx}+\beta \widetilde{w}_{x}=\widetilde{f}_{\varepsilon}(\widetilde{w})
\quad \mbox{for } x\in[0,\widetilde{c}^*_r t],\ t>0,
$$
and
$$
\widetilde{w} (t,\widetilde{c}^*_r  t)=0, \ \ \widetilde{c}^*_r =- \mu \widetilde{w}_x
(t,\widetilde{c}^*_r t) \ \quad  \mbox{for } t\geq0.
$$

Since we are considering the spreading case, we have
$\lim_{t\to\infty}u(t,\cdot)=1$ locally uniformly in $\mathbb{R}$. In particular,
$$
u(t,0)> 1-\varepsilon \quad \mbox{for } t>T
$$
for some $T>0$. Thus $(\widetilde{w} (t,x), \widetilde{c}^*_r t)$ is a lower solution of \eqref{p1}
on $\{(t,x)\mid x\in [0,\widetilde{c}^*_r t],\ t>0\}$ by comparison principle (cf. \cite{DuLin, DuLou}),
and
$$
\widetilde{c}^*_r t \leq h(t+T),\quad \widetilde{w} (t,x) \leq u(t+T,x)
\quad \mbox{ in } \{(t,x)\mid x\in [0,\widetilde{c}^*_r t],\ t>0\}.
$$
This implies that
\begin{equation}\label{est-lower}
\liminf\limits_{t\to \infty}   \frac{h(t)}{t} \geq \widetilde{c}^*_r .
\end{equation}

Next we estimate the upper bound of the rightward spreading speed. Consider the problem
$$
\eta'(t)=f(\eta)\ \ (t>0), \qquad \eta(0)=\|u_0\|_\infty+1.
$$
A simple comparison shows that
$$
u(t,x)\leq\eta(t):= \Big( 1-\frac{\|u_0\|_\infty}{\|u_0\|_\infty+1}e^{-t} \Big)^{-1}
\quad \mbox{ for } x\in [g(t),h(t)],\ t>0.
$$
Hence for any small $\varepsilon>0$, there exists $\widehat{T}>0$ such that
$$
u(t,x)\leq 1+\frac{\varepsilon}{2}\quad \mbox{ for } x\in [0,h(t)],\ t\geq \widehat{T}.
$$
Recall that $(\widehat{c}^*_r, \widehat{q}^*_r (z))$ is a solution of problem
\eqref{prob-q-infty} with $f$ replaced by $\widehat{f}_\varepsilon$ and $\widehat{q}^*_r(\infty)=1+\varepsilon$.
Hence there exists $\widehat{x} > h(\widehat{T})$ large such that
$$
u(\widehat{T}, x) \leq 1+ \frac{\varepsilon}{2} <
\widehat{q}_{r}^{*}(\widehat{x} -x) \quad \mbox{ for } x\in [0,h(\widehat{T})].
$$
Define
$$
\widehat{w} (t,x):=\widehat{q}_{r}^{*}(\widehat{c}_{r}^{*}t +\widehat{x} - x)\quad
\mbox{ for } x\in[0,\widehat{c}^*_r t +\widehat{x}],\ t>0.
$$
Then $(\widehat{w}, \widehat{c}_{r}^{*} t+\widehat{x})$ is an upper solution of \eqref{p1}
on $\{(t,x)\mid x\in [0,h(t+\widehat{T})],\ t>0 \}$, and by the
comparison principle (cf. \cite{DuLin, DuLou}) we have
$$
h(t+\widehat{T})\leq \widehat{c}^*_r t +\widehat{x}, \quad u(t+\widehat{T}, x) \leq \widehat{w} (t,x)
\quad \mbox{ for } x\in [0,h(t+\widehat{T})] \mbox{ and } t>0.
$$
This implies that
\begin{equation}\label{est-upper}
\limsup\limits_{t\to \infty}   \frac{h(t)}{t} \leq \widehat{c}^*_r .
\end{equation}

Since the limits \eqref{est-lower} and \eqref{est-upper} hold for any small $\varepsilon >0$, we have
$$
\lim_{t\to\infty}\frac{h(t)}{t}=c_{r}^{*}
$$
by Proposition \ref{cl-cr}. The leftward spreading speed
$$
\lim_{t\to\infty} \frac{-g(t)}{t}=c_{l}^{*}
$$
is proved similarly.

In \cite{DuLin, DuLou}, the authors considered problem \eqref{p0}, that is, problem
\eqref{p1} {\it without advection term} (i.e. $\beta= 0$). Among others, they showed that
the asymptotic spreading speed  is characterized
by the following problem
\begin{equation}\label{P-beta=0}
\frac{\mathrm{d}p}{\mathrm{d}q}=\frac{c}{d}-\frac{f(q)}{dp}\ \ (q<1),\qquad p(1^-) =0.
\end{equation}
Using a similar phase plane analysis as above, the authors in \cite{DuLou} proved
that problem \eqref{P-beta=0} has a solution $(c, P^c(q))$ for each $c\in [0, 2\sqrt{d})$.
Moreover, they proved that $P^c (0) \searrow 0$ as $c\nearrow c^0 := 2\sqrt{d}$;
\eqref{P-beta=0} has a unique solution $(c^*, P^{c^*} (q))$ such that
$\mu P^{c^*} (0)= c^* $. This $c^*$
is nothing but the rightward and leftward spreading speeds (cf. \cite{DuLin, DuLou}).

Combining the above phase plane analysis we have the following conclusions:
\begin{enumerate}
  \item $P_l^{c-\beta}(0)=P^c(0)=P_r^{c+\beta}(0)$ for all $c\in [\beta, 2\sqrt{d})$;
  \item $P_l^{c}(0)$ (resp. $P^c(0), P_r^{c}(0)$) is continuous and strictly
  decreasing in $c\in [0, 2\sqrt{d}-\beta]$ (resp. $c\in [0,2\sqrt{d}]$, $c\in [0,2\sqrt{d}+\beta]$).
\end{enumerate}
Define three new functions $\gamma_r(c),\ \gamma(c)$ and $\gamma_l(c)$ by
$$
\gamma_r (c) := P^c_r (0)\ \mbox{for } c\in [0, 2\sqrt{d}+\beta),\qquad
\gamma (c) := P^c (0)\ \mbox{for } c\in [0, 2\sqrt{d})
$$
and
$$
\gamma_l (c) := P^c_l (0)\ \mbox{for } c\in [0, 2\sqrt{d}-\beta).
$$
Then, in the $c,\gamma$-plane their graphes lie in the first quadrant (see Figure 1) and these
graphes contact the straight line $\gamma = \frac{c}{\mu}$ at points
$(c^*_l, \frac{c^*_l}{\mu})$, $(c^*,\frac{c^*}{\mu})$ and $(c^*_r,\frac{c^*_r}{\mu})$, respectively.
Therefore, $c^*_l < c^*< c^*_r$. This completes the proof of Theorem \ref{speed}.
\hfill $\Box$

\bigskip
{\it Proof of Theorem \ref{relationship}}. The conclusions in Theorem \ref{relationship}
follow from a simple analysis on the relations among the graphes of $\gamma_r(c), \gamma(c),
\gamma_l(c)$ and $c/\mu$ in Figure 1.\hfill $\Box$

\end{document}